\documentclass[preprint,12pt]{elsarticle}
\usepackage{hyperref}
\usepackage{tabularx}

\usepackage{amsmath}
\usepackage{amssymb}
\usepackage{amsthm}
\newtheorem{thm}{Theorem}[section]
 \newtheorem{cor}[thm]{Corollary}
 \newtheorem{lem}[thm]{Lemma}
 
 \theoremstyle{definition}
 \theoremstyle{rem}
 \newtheorem{rem}{Remark}

  \theoremstyle{exa}

\journal{}

\begin{document}

\begin{frontmatter}



\title{Optimal estimates for mappings admitting general Poisson representations in the unit ball}

\author{Deguang Zhong\fnref{label2}}
\fntext[label2]{huachengzhon@163.com}
\address{1. Institute of Applied Mathematics, Shenzhen Polytechnic University,\\
Shenzhen, Guangdong, 518055, P. R. China}

\author{Fangming Cai\fnref{label3}}
\fntext[label3]{caifangming@hnu.edu.cn}
\address{2. School of Mathematics, Hunan University,\\
Changsha, Hnan, 410082, P. R. China}

\author{Dongping Wei\corref{cor1}\fnref{label4}}
\cortext[cor1]{Corresponding author.}
\fntext[label4]{dongpingwei2023@163.com}
\address{3. Institute of Applied Mathematics, Shenzhen Polytechnic University,\\
Shenzhen, Guangdong, 518055, P. R. China}


\begin{abstract}
  Suppose that $1<p\leq\infty$ and  $\varphi\in L^{p}(\mathbb{B}^{n},\mathbb{R}^{n}).$  In this note, we use H\"{o}lder inequality and some basic properties  of hypergeometric functions to establish the sharp constant $C_{p}$ and function $C_{p}(x)$ in the following inequalities
$$|u(x)|\leq \frac{C_{p}}{(1-|x|^{2})^{(n-1)/p}}\cdot||\varphi||_{L^{p}}$$
  and
  $$|u(x)|\leq \frac{C_{p}(x)}{(1-|x|^{2})^{(n-1)/p}}\cdot||\varphi||_{L^{p}},$$
  where  $u$ are those mapping from the unit ball $\mathbb{B}^{n}$ into $\mathbb{R}^{n}$ admitting general Poisson representations. The obtained results generalize and extend some known results from harmonic mappings (\cite[Proposition 6.16]{ABR92} and \cite[Theorems 1.1 and 1.2]{DM12}) and hyperbolic harmonic mappings (\cite[Theorems 1.1 and 1.2]{CJLK20}).
\end{abstract}

\begin{keyword}
Harmonic mapping\sep hyperbolic harmonic mapping\sep general Poisson representations
 \MSC[2020] 31B05,31B10, 42B30
\end{keyword}

\end{frontmatter}


\section{Introduction and main results}\label{sec1}
 Let $\mathbb{B}^{n}$ be the unit ball of $\mathbb{R}^{n}$ and $\mathbb{S}^{n-1}$ be its boundary. For any $x=(x_{1}, x_{2},\ldots, x_{n})\in\mathbb{R}^{n},$ its vector norm is defined by $|x|=(\sum_{i=1}^{n}x_{i}^{2})^{1/2}.$  For $a, b, c\in\mathbb{R}, c\neq0, -1, -2,\cdot\cdot\cdot,$ the hypergeometric function is defined by
$$F(a, b; c; x)=\sum_{k=0}^{\infty}\frac{(a)_{k}\cdot(b)_{k}}{(c)_{k}}\cdot\frac{x^{k}}{k!},\;|x|<1,$$
where $(a)_{k}=\Gamma(a+k)/\Gamma(a)$ and $\Gamma$ is the Gamma function; cf. \cite{AAR99}. The following two elemental properties on hypergeometric function is well known; cf. \cite[Chapter 2]{AAR99}.
\begin{equation}\label{euqa2}
F(a, b; c; x)=(1-x)^{c-a-b}\cdot F(c-a, c-b; c; x).
\end{equation}
\begin{equation}\label{euqa3}
\lim _{x\rightarrow1^{-}}F(a, b; c; x)=\frac{\Gamma(c)\Gamma(c-a-b)}{\Gamma(c-a)\Gamma(c-b)}\;\;{\rm if}\; c-a-b>0.
\end{equation}
For $\alpha,\beta\in\mathbb{R},\beta>0.$ Then the general Poisson kernel \cite{KMM21} is defined by
\begin{equation}
P_{\alpha,\beta}(x,\eta)=\frac{(1-|x|^{2})^{\alpha}}{|x-\eta|^{\beta}},\;\;x\in\mathbb{B}^{n} \;{\rm and}\;\;\eta\in\mathbb{S}^{n-1}.
\end{equation}
In the following discussion, we always consider the positive integers $n$ which greater than or equal to $3$.  In this paper, we are interesting in those mappings admitting the following general Poisson representations:
 \begin{equation}\label{equation0}
u_{\alpha,\beta}[\varphi](x)=c_{n,\beta}\cdot\int_{\mathbb{S}^{n-1}}P_{\alpha,\beta}(x,\eta)\varphi(\eta)d\sigma(\eta),
\end{equation}
where $\varphi\in L^{1}(\mathbb{S}^{n-1},\mathbb{R}^{n}),$ $c_{n,\beta}:=\frac{\Gamma(\frac{\beta}{2})\Gamma(\frac{\beta}{2}-\frac{n}{2}+1)}{\Gamma(\frac{n}{2})\Gamma(\beta-n+1)}$ and $\sigma$ is the surface measure on $\mathbb{S}^{n-1}$ normalized by $\sigma(\mathbb{S}^{n-1})=1.$  It is noted that  many classical functions defined in the unit ball can be written as the form in (\ref{equation0}). For example, a harmonic mapping $u$ defined in the unit ball can be represented as $u(x)=u_{1,n}[\varphi](x),$ where $\varphi\in L^{1}(\mathbb{S}^{n-1}, \mathbb{R}^{n}).$
If $\alpha=n-1, \beta=2(n-1),$ then the function $u(x)=u_{n-1,2(n-1)}[\varphi](x)$ is hyperbolic harmonic mapping. It was show by Liu and Peng \cite{LP04} that the following Dirichlet problem
\begin{equation}\label{euq1}
\left\{
\begin{aligned}
\Delta_{\gamma}u&=0\;\;\;\;{\rm in}\;\;\;\mathbb{B}^{n},\\
u&=f\;\;\;\;{\rm on}\;\;\mathbb{S}^{n-1}
\end{aligned}
\right.
\end{equation}
has a solution for all $f\in \mathcal{C}(\mathbb{S}^{n-1})$ if and only if $\gamma>-1/2,$ where
$$\Delta_{\gamma}=(1-|x|^{2})\cdot\left[\frac{1-|x|^{2}}{4}\cdot\sum_{i}\frac{\partial^{2}}{\partial x_{i}^{2}}+\gamma\sum_{i}x_{i}\cdot\frac{\partial}{\partial x_{i}}+\gamma(\frac{n}{2}-1-\gamma)\right].$$
In this case the solution is unique and is represented  as
$u(x)=u_{1+2\gamma,n+2\gamma}[f](x).$

In this paper, our aim is to obtain the sharp constants $C_{p}$ and sharp function $C_{p}(x)$ in the following inequalities
 \begin{equation}\label{ine12}
|u_{\alpha,\beta}[\varphi](x)|\leq \frac{C_{p}(x)}{(1-|x|^{2})^{(n-1)/p}}\cdot||\varphi||_{L^{p}}
\end{equation}
 and
  \begin{equation}\label{ine12}
|u_{\alpha,\beta}[\varphi](x)|\leq \frac{C_{p}}{(1-|x|^{2})^{(n-1)/p}}\cdot||\varphi||_{L^{p}},
\end{equation}
 where $\varphi\in L^{p}(\mathbb{S}^{n-1},\mathbb{R}^{n}).$

A similar result proved in \cite[Lemma 5.1.1]{MP04}, which states that if $f\in H^{p},0<p\leq\infty,$ then $$|f(z)|\leq(1-|z|^{2})^{-1/p}||f||_{p}.$$ Here $H^{p}$ stands for the Hardy space consisting of analytic functions $f$ and satisfying
$$||f||_{p}:=\sup_{r<1}\left(\frac{1}{2\pi}\int_{-\pi}^{\pi}|f(re^{i\theta})|^{p}\right)^{1/p}<\infty.$$
Suppose $1\leq p \leq\infty.$ Let $h^{p}(\mathbb{B}^{n})$ be the harmonic Hardy spaces on the unit ball $\mathbb{B}^{n},$ which the function $u\in h^{p}(\mathbb{B}^{n})$ satisfying the following conditions
$$||u||_{p}^{p}:=||u||_{h^{p}(\mathbb{B}^{n})}^{p}=\sup_{0<r<1}\int_{\mathbb{S}^{n-1}}|u(r\zeta)|^{p}d\sigma(\zeta)<\infty.$$
For the case of $p=q=2,$ the following sharp estimate
$$|u(x)|\leq\sqrt{\frac{1+|x|^{2}}{(1-|x|^{2})^{n-1}}}\cdot||u||_{h^{2}(\mathbb{B}^{n})}$$
was obtained in \cite[Proposition 6.23]{ABR92}. In \cite{DM12},  Kalaj and Markovi\'{c} obtained  the following optimal estimates for harmonic functions $u_{1,n}[\varphi]$ in the unit ball.
\begin{thm}{\rm\cite[Theorems 1.1 and 1.2]{DM12}}\label{thm1}
 Let $1<p\leq\infty$ and $q$ be its conjugate. For all $u\in h^{p}(\mathbb{B}^{n})$ and $x\in \mathbb{B}^{n},$ we have the following sharp inequalities
\begin{equation}\label{ine12}
|u(x)|\leq \frac{C_{p}(x)}{(1-|x|^{2})^{(n-1)/p}}\cdot||u||_{h^{p}(\mathbb{B}^{n})}
\end{equation}
and
\begin{equation}\label{ine12}
|u(x)|\leq \frac{C_{p}}{(1-|x|^{2})^{(n-1)/p}}\cdot||u||_{h^{p}(\mathbb{B}^{n})},
\end{equation}
 where
 \begin{equation}\label{ine88}
 C_{p}(x)= \left(F\left(\frac{n-nq}{2}, -1+n-\frac{nq}{2}; \frac{n}{2};|x|^{2}\right)\right)^{1/q}
 \end{equation}
 and
 \begin{equation}\label{ine888}
  C_{p}=\left\{
\begin{aligned}
&1,       &      & {\rm if}\; {q\leq\frac{2(n-1)}{n}},\\
&\left(\frac{2^{nq-n}\Gamma(\frac{n}{2})\Gamma(\frac{nq-n+1}{2})}{\sqrt{\pi}\Gamma(\frac{nq}{2})}\right)^{1/q},       &      &{\rm if}\; {q>\frac{2(n-1)}{n}.}
\end{aligned} \right.
\end{equation}
\end{thm}

Later, Chen and Kalaj \cite{CJLK20} derived the following optimal estimates for hyperbolic Poisson integrals of functions $u_{n-1,2n-2}[\varphi]$ in the unit ball.
\begin{thm}{\rm\cite[Theorems 1.1 and 1.2]{CJLK20}}\label{thm2}
 Let $1<p\leq\infty$ and $q$ be its conjugate. If $u=u_{n-1, 2n-2}[\varphi]$ and $\varphi\in L^{p}(\mathbb{S}^{n-1},\mathbb{R}^{n}),$ then for any $x\in \mathbb{B}^{n},$ we have the following sharp inequalities
\begin{equation}\label{ine12}
|u(x)|\leq \frac{C_{p}(x)}{(1-|x|^{2})^{(n-1)/p}}\cdot||\varphi||_{L^{p}}
\end{equation}
and
\begin{equation}\label{ine12}
|u(x)|\leq \frac{C_{p}}{(1-|x|^{2})^{(n-1)/p}}\cdot||\varphi||_{L^{p}},
\end{equation}
 where
 \begin{equation}\label{ine99}
 C_{p}(x)= \left(F\left(-(n-1)(q-1), \frac{n}{2}+q-nq; \frac{n}{2};|x|^{2}\right)\right)^{1/q}
 \end{equation}
 and
  \begin{equation}\label{ine999}
 C_{p}= \left(\frac{\Gamma(\frac{n}{2})\Gamma((2q-1)(n-1))}{\Gamma(\frac{n}{2}+(q-1)(n-1))\Gamma(q(n-1))}\right)^{1/q}.
 \end{equation}
\end{thm}

 In this paper, we will povide a simple and direct  method to established the following optimal estimates for mappings admitting general Poisson representations $u_{\alpha,\beta}[\varphi]$  in the unit ball. 
\begin{thm}\label{thm3}
 Let $1<p\leq\infty$ and $q$ be its conjugate. Suppose that $\beta-\alpha\geq n-1,\alpha\geq1.$ If $u(x)=(1-|x|^{2})^{\beta-\alpha-\frac{n-1}{q}}\cdot u_{\alpha,\beta}[\varphi](x)$ and $\varphi\in L^{p}(\mathbb{S}^{n-1},\mathbb{R}^{n}).$ Then for any $x\in \mathbb{B}^{n},$ we have the following sharp inequality
\begin{equation}\label{ine12}
|u(x)|\leq c_{n,\beta} \cdot \left(F\left(\frac{n-q\beta}{2}, n-1-\frac{q\beta}{2};\frac{n}{2}; |x|^{2}\right)\right)^{1/q}\cdot||\varphi||_{L^{p}}.
\end{equation}
\end{thm}

\begin{thm}\label{thm4}
 Let $1<p\leq\infty$ and $q$ be its conjugate. Suppose that $\beta-\alpha\geq n-1,\alpha\geq1.$ If $u(x)=(1-|x|^{2})^{\beta-\alpha-\frac{n-1}{q}}\cdot u_{\alpha,\beta}[\varphi](x)$  and $\varphi\in L^{p}(\mathbb{S}^{n-1},\mathbb{R}^{n}).$ Then for any $x\in \mathbb{B}^{n},$ we have the following sharp inequality
\begin{equation}\label{ine19}
|u(x)|\leq C_{p}\cdot||\varphi||_{L^{p}}.
\end{equation}
 In the case that $n\leq\beta<2(n-1),$ we have
\begin{equation}\label{ine111111}
  C_{p}=\left\{
\begin{aligned}
&c_{n,\beta},      &   &    {\rm if}\;{q\leq\frac{2(n-1)}{\beta}},\\
&c_{n,\beta}\cdot\left(\frac{\Gamma(\frac{n}{2})\Gamma(q\beta -n+1)}{\Gamma(\frac{q\beta}{2})\Gamma(\frac{\beta q-n+2}{2})}\right)^{1/q},       &      &{\rm if}\; {q>\frac{2(n-1)}{\beta}.}
\end{aligned} \right.
\end{equation}
 In the case that $\beta\geq2(n-1),$ we have
\begin{equation}
C_{p}=c_{n,\beta}\cdot\left(\frac{\Gamma(\frac{n}{2})\Gamma(q\beta -n+1)}{\Gamma(\frac{q\beta}{2})\Gamma(\frac{\beta q-n+2}{2})}\right)^{1/q}.
\end{equation}
\end{thm}

\begin{rem}\label{rem1}
  If  $p=\infty,$ then the proofs of Theorems \ref{thm3} and \ref{thm4}  are trivial. Here, we leave the readers to check for this case.
\end{rem}
\begin{rem}
  If let $\alpha=1, \beta=n$ in {\rm(\ref{ine12})} and  {\rm(\ref{ine111111})}, that we can derive equations {\rm(\ref{ine88})} and  {\rm(\ref{ine888})}. The equations in {\rm(\ref{ine99})} and  {\rm(\ref{ine999})} can be obtained be letting $\alpha=n-1, \beta=2(n-1)$ in  equations {\rm(\ref{ine12})} and  {\rm(\ref{ine111111})}. If we consider those mappings satisfying Dirichlet problem {\rm(\ref{euq1})}, then we get the following corollary.
\end{rem}
\begin{cor}
  Suppose that  $f\in \mathcal{C}(\mathbb{S}^{n-1})\bigcap L^{p}(\mathbb{S}^{n-1},\mathbb{R}^{n}).$ Let $u$ be a function from unit ball satisfying the Dirichlet problem {\rm(\ref{euq1})}.
Then for any $x\in \mathbb{B}^{n},$ we have the following inequalities
\begin{equation}\label{ine77}
|u(x)|\leq \frac{C_{p}(x)}{(1-|x|^{2})^{(n-1)/p}}\cdot||\varphi||_{L^{p}}
\end{equation}
 and
\begin{equation}\label{ine120}
|u(x)|\leq \frac{C_{p}}{(1-|x|^{2})^{(n-1)/p}}\cdot||\varphi||_{L^{p}},
\end{equation}
where
 \begin{equation}\label{ine1777}
 C_{p}(x)=c_{n,n+2\gamma} \cdot F\left(\frac{(1-q)n}{2}-q\gamma, \frac{(2-q)n}{2}-q\gamma-1;\frac{n}{2}; |x|^{2}\right)^{1/q}
 \end{equation}
and in the case that $0\leq\gamma<\frac{n}{2}-1,$
\begin{equation}\label{ine11111111111}
  C_{p}=\left\{
\begin{aligned}
&c_{n,n+2\gamma},      &   &    {\rm if}\;{q\leq\frac{2(n-1)}{n+2\gamma}},\\
&c_{n,n+2\gamma}\cdot\left(\frac{\Gamma(\frac{n}{2})\Gamma((q-1)n+2q\gamma+1)}{\Gamma(\frac{nq+2\gamma q}{2})
\Gamma\left(\frac{(q-1)n}{2}+q\gamma+1\right)}\right)^{1/q},       &      &{\rm if}\; {q>\frac{2(n-1)}{n+2\gamma}}
\end{aligned} \right.
\end{equation}
and in the case that $\gamma\geq\frac{n}{2}-1,$
\begin{equation}
C_{p}=c_{n,n+2\gamma}\cdot\left(\frac{\Gamma(\frac{n}{2})\Gamma((q-1)n+2q\gamma+1)}{\Gamma(\frac{nq+2\gamma q}{2})
\Gamma\left(\frac{(q-1)n}{2}+q\gamma+1\right)}\right)^{1/q}.
\end{equation}
\end{cor}
The rest of the paper is organized as follows: In Sect. \ref{sec2}, we will make some preparations which will be used in proving our main results. In Sect. \ref{sec3}, the proof of Theorem \ref{thm3} is given.  The last Section will be devoted to proving Theorem \ref{thm4}.

\section{Preliminaries}\label{sec2}
The following lemma concerning the monotonicity of hypergeometric functions was proved in \cite{O14}.
\begin{lem}{\rm \cite[Lemma 1.2]{O14}}\label{lem1}
Suppose that $c>0, a\leq c, b\leq c$ and $ab\leq0\;(ab\geq0).$ Then the hypergeometric function $F(a, b; c;\cdot)$ is decreasing {\rm(}increasing{\rm)} on $(0, 1).$
\end{lem}

The proof of following lemma can be found in \cite[Lemma 2.1]{LP04}.
\begin{lem}{\rm\cite[Lemma 2.1]{LP04}}\label{lem09}
For  $x\in\mathbb{B}^{n}$ and $\lambda\in \mathbb{C}.$ Then we have
\begin{equation}\label{ine20}
\int_{\mathbb{S}^{n-1}}\frac{1}{|x-\eta|^{2\lambda}}d\sigma(\eta)=F(\lambda, \lambda-\frac{n}{2}+1; \frac{n}{2}; |x|^{2}).
\end{equation}
\end{lem}

\section{The proof of Theorem \ref{thm3}}\label{sec3}
\begin{proof}
 By Remark \ref{rem1}, we here only consider the case when $1<p<\infty$. By the classical H\"{o}lder's inequality, Lemma \ref{lem09} and equation (\ref{euqa2}), we get
 \begin{equation}\label{ontvnium}
\begin{aligned}
&|u(x)|=c_{n,\beta}\cdot(1-|x|^{2})^{\beta-\alpha-\frac{n-1}{q}}\Big|\int_{\mathbb{S}^{n-1}}P_{\alpha,\beta}(x,\eta)\varphi(\eta)d\sigma(\eta)\Big|\\
\leq&c_{n,\beta}\cdot\left(\int_{\mathbf{S}^{n-1}}\left(\frac{(1-|x|^{2})^{\beta-\frac{n-1}{q}}}{|x-\eta|^{\beta}}\right)^{q}d\sigma(\eta)\right)^{1/q}\cdot \left(\int_{\mathbf{S}^{n-1}}|\varphi(\eta)|^{p}d\sigma(\eta)\right)^{1/p}\\
=&c_{n,\beta}\cdot(1-|x|^{2})^{\beta-\frac{n-1}{q}}\cdot\left(\int_{\mathbf{S}^{n-1}}\frac{1}{|x-\eta|^{q\beta}}d\sigma(\eta)\right)^{1/q}\cdot||\varphi||_{L^{p}}\\
=&c_{n,\beta}\cdot(1-|x|^{2})^{\beta-\frac{n-1}{q}}\cdot \left(F\left(\frac{q\beta}{2}, \frac{q\beta}{2}-\frac{n}{2}+1;\frac{n}{2}; |x|^{2}\right)\right)^{1/q}\cdot||\varphi||_{L^{p}}\\
=&c_{n,\beta}\cdot(1-|x|^{2})^{\beta-\frac{n-1}{q}}\cdot(1-|x|^{2})^{\frac{n-q\beta-1}{q}}\cdot \left(F\left(\frac{n-q\beta}{2}, n-1-\frac{q\beta}{2};\frac{n}{2}; |x|^{2}\right)\right)^{1/q}\cdot||\varphi||_{L^{p}}\\
=&c_{n,\beta}\cdot \left(F\left(\frac{n-q\beta}{2}, n-1-\frac{q\beta}{2};\frac{n}{2}; |x|^{2}\right)\right)^{1/q}\cdot||\varphi||_{L^{p}}.\\
\end{aligned}
\end{equation}
In order to prove (\ref{ontvnium}) is sharp, we take $$\varphi_{0}(\eta)=\left(0,0,\ldots,\left[\frac{(1-|x|^{2})^{\beta-\frac{n-1}{q}}}{|x-\eta|^{\beta}}\right]^{q/p}\right),$$ where $x\in \mathbb{B}^{n},\eta\in \mathbb{S}^{n-1}.$ 
We first verify that $\varphi_{0}\in L^{p}(\mathbb{S}^{n-1},\mathbb{R}^{n}).$ This is because 
\begin{equation}
\begin{aligned}
(||\varphi_{0}||_{L^{p}}(\mathbb{S}^{n-1}))^{p}=&\int_{\mathbb{S}^{n-1}}\left[\frac{(1-|x|^{2})^{\beta-\frac{n-1}{q}}}{|x-\eta|^{\beta}}\right]^{q} d\sigma(\eta)\\
=&(1-|x|^{2})^{q\beta-(n-1)}\int_{\mathbb{S}^{n-1}}\frac{1}{|x-\eta|^{q\beta}}d\sigma(\eta)\\
=&(1-|x|^{2})^{q\beta-(n-1)}\cdot F\left(\frac{q\beta}{2}, \frac{q\beta}{2}-\frac{n}{2}+1;\frac{n}{2}; |x|^{2}\right)\\
=&F\left(\frac{n-q\beta}{2}, n-1-\frac{q\beta}{2};\frac{n}{2}; |x|^{2}\right)\leq C_{p}<\infty,\\
\end{aligned}
\end{equation}
where $C_{p}$ is given by (\ref{ine19}) and (\ref{ine111111}).

 Since
\begin{equation}
\begin{aligned}
|u(x)|&=c_{n,\beta}\cdot\int_{\mathbb{S}^{n-1}}\frac{(1-|x|^{2})^{\beta-\frac{n-1}{q}}}{|x-\eta|^{\beta}}\cdot \left[\frac{(1-|x|^{2})^{\beta-\frac{n-1}{q}}}{|x-\eta|^{\beta}}\right]^{q/p} d\sigma(\eta)\\
&=c_{n,\beta}\cdot\int_{\mathbb{S}^{n-1}}\left[\frac{(1-|x|^{2})^{\beta-\frac{n-1}{q}}}{|x-\eta|^{\beta}}\right]^{q} d\sigma(\eta)\\
\end{aligned}
\end{equation}
and
\begin{equation}
\begin{aligned}
&c_{n,\beta}\cdot\left(\int_{\mathbf{S}^{n-1}}\left[\frac{(1-|x|^{2})^{\beta-\frac{n-1}{q}}}{|x-\eta|^{\beta}}\right]^{q}d\sigma(\eta)\right)^{1/q}\cdot \left(\int_{\mathbf{S}^{n-1}}|\varphi_{0}(\eta)|^{p}d\sigma(\eta)\right)^{1/p}\\
=&c_{n,\beta}\cdot\left(\int_{\mathbf{S}^{n-1}}\left[\frac{(1-|x|^{2})^{\beta-\frac{n-1}{q}}}{|x-\eta|^{\beta}}\right]^{q}d\sigma(\eta)\right)^{1/q}\cdot \left(\int_{\mathbf{S}^{n-1}}\left[\frac{(1-|x|^{2})^{\beta-\frac{n-1}{q}}}{|x-\eta|^{\beta}}\right]^{q}d\sigma(\eta)\right)^{1/p}\\
=&c_{n,\beta}\cdot\int_{\mathbb{S}^{n-1}}\left[\frac{(1-|x|^{2})^{\beta-\frac{n-1}{q}}}{|x-\eta|^{\beta}}\right]^{q} d\sigma(\eta),\\
\end{aligned}
\end{equation}
we get $$|u(x)|=c_{n,\beta}\cdot\left(\int_{\mathbf{S}^{n-1}}\left(\frac{(1-|x|^{2})^{\beta-\frac{n-1}{q}}}{|x-\eta|^{\beta}}\right)^{q}d\sigma(\eta)\right)^{1/q}\cdot||\varphi_{0}||_{L^{p}}.$$
This shows  (\ref{ine12}) is sharp. Hence, the proof is finished.
\end{proof}

\section{The proof of Theorem \ref{thm4}}\label{sec4}
\begin{proof}
  Let 
  \begin{equation}\label{sjiknscsccs}
   C_{p}(x):=c_{n,\beta} \cdot \left(F\left(\frac{n-q\beta}{2}, n-1-\frac{q\beta}{2};\frac{n}{2}; |x|^{2}\right)\right)^{1/q}.
  \end{equation}
  Here we also only consider the case when $1<p<\infty$. The main ideal of proof of Theorem \ref{sec4} is to study the monotonicity of function $$\psi(r):=F\left(\frac{n-q\beta}{2}, n-1-\frac{q\beta}{2};\frac{n}{2}; r\right),r\in[0,1).$$ This can be done by virtue of Lemma \ref{lem1}. To do this, it is noted that  $\frac{n}{2}>\frac{n-q\beta}{2}, n/2>n-1-\frac{q\beta}{2}.$  Next, we will discuss in the following two cases:

Case (i): If $n\leq\beta<2(n-1).$ Since in this case, we have
 \begin{equation}\label{ine1111}
  \left\{
\begin{aligned}
&(\frac{n-q\beta}{2})(n-1-\frac{q\beta}{2})\leq0,       &      & {\rm if}\; {q\leq\frac{2(n-1)}{\beta}},\\
&(\frac{n-q\beta}{2})(n-1-\frac{q\beta}{2})>0,       &      &{\rm if}\; {q>\frac{2(n-1)}{\beta}.}
\end{aligned} \right.
\end{equation}
Hence, we see from Lemma \ref{lem1} that $\psi$ is monotonically decreasing (increasing) on $[0,1)$ if $q\leq 2(n-1)/\beta$ ($q>2(n-1)/\beta$).   So there holds
 \begin{equation}\label{ine3333}
  C_{p}=\max_{0\leq|x|<1}C_{p}(x)=\left\{
\begin{aligned}
&C_{p}(0),       &      & {\rm if}\; {q\leq\frac{2(n-1)}{\beta}},\\
&C_{p}(1),       &      &{\rm if}\; {q>\frac{2(n-1)}{\beta}.}
\end{aligned} \right.
\end{equation}

 Case (ii): If $2(n-1)\leq\beta.$ In this case, we have $$(\frac{n-q\beta}{2})(n-1-\frac{q\beta}{2})>0\;\;\;{\rm for \;\;all}\;\;\; q>1,$$
 which implies that $\psi$ is monotonically increasing on $[0,1)$ for all  $q>1$ according to Lemma \ref{lem1}. Hence, in this case, we have $C_{p}=C_{p}(1).$

At last, we will get the values of $C_{p}(0)$ and $C_{p}(1).$ It is obvious that $C_{p}(0)=c_{n,\beta}.$ In addition, by combining equation (\ref{euqa3}) with equation (\ref{sjiknscsccs}), we derive
\begin{equation}\label{ine21}
  C_{p}(1)=c_{n,\beta}\cdot\left(\frac{\Gamma(\frac{n}{2})\Gamma(q\beta -n+1)}{\Gamma(\frac{q\beta}{2})\Gamma(\frac{\beta q-n+2}{2})}\right)^{1/q}.
\end{equation}
This finished the proof of Theorem \ref{thm4}.
\end{proof}

\paragraph{Acknowledgements}   The first author was supported by Supporting Foundation of Shenzhen Polytechnic University (No. ZX2023000301)
and Research Foundation of Shenzhen Polytechnic University (No. 6023312032K). The third author was supported by Guangdong Province Higher Vocational Education Teaching Reform Research and Practice Project of  in 2020 (No. JGGZKZ2020167).











\end{document}